\documentclass[12pt]{article}
\usepackage{amsmath}
\newcommand{\be}{\begin{equation}}
\newcommand{\ee}{\end{equation}}
\newcommand{\ba}{\begin{eqnarray}}
\newcommand{\ea}{\end{eqnarray}}
\newcommand{\ban}{\begin{eqnarray*}}
\newcommand{\ean}{\end{eqnarray*}}

 \newcommand{\qed}{\hspace*{\fill}\rule{3mm}{3mm}\quad}
\newcommand{\Pf}{\noindent  {\em Proof.} }

\newcommand{\vol}{\mathrm{Vol}}

\newcommand{\sect}[1]{\section{#1}  \setcounter{equation}{0}}

\newtheorem{lem}{Lemma}[section]
 
\begin{document}
\newtheorem{defn}[lem]{Definition}
\newtheorem{theo}[lem]{Theorem}
\newtheorem{cor}[lem]{Corollary}
\newtheorem{prop}[lem]{Proposition}
\newtheorem{rk}[lem]{Remark}
\newtheorem{ex}[lem]{Example}
\newtheorem{note}[lem]{Note}
\newtheorem{conj}[lem]{Conjecture}

\title{A Neumann Type Maximum Principle for the Laplace Operator on Compact  
Riemannian Manifolds}
\author{Guofang Wei and Rugang Ye \\ {\small Department  of Mathematics} \\
{\small University of California, Santa  Barbara}}
\date{}
\maketitle
 \begin{abstract}
 In this paper we present a proof of a Neumann type
maximum principle for the Laplace operator on compact Riemannian
manifolds. A key point is the simple geometric nature of the constant in
the a priori estimate
of this maximum principle.  In particular, this maximum principle can  be
applied to manifolds 
with Ricci curvature bounded from below and diameter 
bounded from above to yield a maximum estimate without dependence on a positive 
lower bound for the volume.
\end{abstract}

\sect{Introduction}
 
The main purpose of this paper is to present a proof of a Neumann type  
maximum principle for 
the Laplace operator on a closed Riemannian manifold. As a key feature of this maximum principle,  
the constant in the maximum estimate depends on the Riemannian manifold only in terms of 
the dimension and the volume-normalized Neumann isoperimetric constant.  This allows us to apply it to 
manifolds with Ricci curvature bounded from  below and diameter bounded from above to obtain a maximum 
principle 
without dependence  on a positive lower bound for the volume.     
 A special case of this 
maximum principle,
namely Theorem C with $\Phi=0$  has been believed to be true and used in 
[P] for establishing an eigenvalue pinching theorem for manifolds with positive Ricci 
curvature.  (The 
accounts in [P] also suggest a belief in a general version.)
 But we cannot 
find any other reference for this maximum principle (the special case or the general case) 
in the literature. A corresponding maximum principle (in various 
formulations) 
for  the Dirichlet 
boundary value problem on a domain  is well-known. But its 
usual proof, which is an  application of Moser iteration based on the Sobolev 
inequality, 
is not suitable for  the Neumann type  problem of this paper for a number of reasons.  In particular, 
the key independence 
from volume lower bound mentioned above requires new
arguments  for our Neumann type problem. 
Another obvious difference is that 
no average of the subsolution appears in the maximum principle for the Dirichlet 
boundary value problem, in contrast to the situation of this paper.

  Consider a closed Riemannian manifold 
$(M, g)$ of dimension $n$,  where $g$ denotes the metric. Let $L^p(M)$ denote 
the $L^p$ space 
of  functions on $M$, $L^p(TM)$ the $L^p$ space of vector fields on $M$, and  
$W^{k, p}(M)$ the $W^{k,p}$ Sobolev space of functions on $M$.  The  $L^p$ 
norm with respect 
to $g$ will be denoted by $\| \cdot \|_p$, i.e.  
\ba
\|f\|_p =\left(\int_M |f|^p\right)^{\frac{1}{p}},  \|\Phi\|_p=\left(\int_M 
|\Phi|^p\right)^{\frac{1}{p}}
\ea
for $f \in  L^p(M)$ and $\Phi \in L^p(TM)$.
(The notation of the volume form of $g$ is  often omitted in this paper.) 
The  following volume-normalized 
$L^p$ norm $ \| \cdot  \|_p^*$ will play an important role in this paper:
\ba
\|f\|_p^* = \left(\frac{1}{vol_g(M)} \int_M  |f|^p\right)^{\frac{1}{p}}, 
\|\Phi\|^*_p  &=&\left(\frac{1}{vol_g(M)}\int_M|\Phi|^p\right)^{\frac{1}{p}},
\ea
where  $vol_g(M)$ denotes the volume of $(M,g)$. 
 
The average of a function $u \in L^1(M)$ on $M$ will be denoted by $u_M$,  
i.e.
\ba
u_M=\frac{1}{vol_g(M)}\int_M u.
\ea
 
For a function $u$ on $M$ we denote its positive part by $u^+$ and its  
negative part by 
$u^-$, i.e. $u^+=\max\{u, 0\}$ and $u^-=\min\{u, 0\}$. 
The Laplace operator 
$\Delta$ is  the negative Laplacian, i.e. $\Delta u=\mbox{div} \nabla u$. Let $C^*_{I,N}(M,g)$ denote 
the volume-normalized Neumann isoperimetric constant, which is defined 
in terms of the Neumann isoperimentric constant $C_{I,N}(M,g)$,  see Section 2.  \\

\noindent {\bf Theorem A} {\it Assume $n \ge 3$. Let  $u$ be a  function in 
$W^{1,\alpha}(M)$ with $\alpha>n$, which satisfies
\ba  \label{sub}
\Delta u \ge f+\mbox{ div }\Phi
\ea
in the weak sense   for a measurable function $f$ on $M$ such that $f^- \in 
L^p(M)$  and 
a  vector field
$\Phi \in L^{2p}(TM)$ with $p>\frac{n}{2}$, i.e. 
\ba  \label{nabla3}
\int_M \nabla u \cdot \nabla \phi \le 
-\int_M f \phi +  \int_M \Phi \cdot \nabla \phi
\ea
for all nonnegative $\phi \in  W^{1,2}(M)$ (equivalently, all nonnegative 
$\phi \in  W^{1,\frac{\alpha}{\alpha-1}}(M)$).     Then we have 
\ba  \label{maxA}
\sup_M \, u \le u_M+ C(n, p, C^*_{N,I}(M,g))  
(\|f^-\|^*_p+ \|\Phi\|^*_{2p})
\ea
with a positive  constant $C(n, p, C^*_{N,I}(M,g))$ depending only on 
$n, p$ and  $C^*_{N,I}(M,g)$.
This constant depends continuously and  
increasingly on $C^*_{N,I}(M,g)$. } \\

The classical strong maximum principle says that $u\equiv u_M$ if $\Delta u  
\ge 0$ (in the weak sense). Theorem A includes this as a special  corollary.  
But the main 
point of Theorem A lies in the quantitative  estimate (\ref{maxA}) and  
the simple geometric nature of the constant  $C(n,p,C^*_{N,I}(M,g))$ in the 
estimate. As emphasized above,  no other 
data from the metric  $g$ such as the volume are involved in this constant. 

In contrast to traditional estimates of the maximum principle type, the estimate 
(\ref{maxA}) is not scaling invariant. In other words, the estimate obtained with respect to 
a rescaled metric and the corresponding rescaled $f$ and $\Phi$ differs from the original 
estimate. This non-invariace is brought into the estimate by a construction in the proof of 
Lemma 4.2, see Remark 3 in Section 4. (For a discussion of the behavior of the estimate 
under rescaling of the metric see Remark 4.)   Without breaking the scaling invariance it would be impossible to obtain 
a maximum estimate in which the constant depends solely on the dimension $n$, the exponent $p$  and the 
volume-normalized Neumann  isoperimetric constant. This is one of the key features of our 
arguments.  (Scaling invariant maximum estimates can also be 
derived, see [Y].)

As a consequence of Theorem A and S.~Gallot's estimate of the volume-normalized  Neumann 
isoperimetric 
constant in [Ga1] (see Theorem \ref{gallot}) we obtain  the following result 
which involves
a lower bound for the Ricci curvature and  an upper bound for the diameter. 
For convenience, we define the {\it  diameter rescaled Ricci curvature} of a 
unit tangent 
vector $v$ to be $\hat  Ric(v,v)=diam_g(M)^2 Ric(v,v)$, where $diam_g(M)$ 
denotes 
the diameter of  $(M,g)$. We 
set  $\kappa_{\hat Ric}=\min_{v\in TM, |v|=1} \hat  Ric(v,v)$ \
and $\hat \kappa_{\hat Ric}=|\kappa_{\hat  Ric}^-|=|\min\{\kappa_{\hat Ric}, 
0\}|$. \\

\noindent {\bf Theorem B}  {\it Assume $n \ge 3$. Let  $u$ be  a function in 
$W^{1,\alpha}(M)$ with $\alpha>n$ satisfying 
\ba
\Delta  u \ge f+\mbox{ div }\Phi
\ea
in the weak sense for a measurable function  $f$ such that $f^- \in L^p(M)$  
and 
a vector field $\Phi \in  L^{2p}(TM)$ with $p>\frac{n}{2}$.   Then we have 
\ba  \label{maxB}
\sup_M \, u \le u_M+ C(n, p, \hat \kappa_{\hat Ric}, diam_g(M))  
(\|f^-\|^*_p+ 
\|\Phi\|^*_{2p}).
\ea
with a positive  constant $C(n, p, \hat \kappa_{\hat Ric}, diam_g(M))$ depending only on 
$n, p, \hat \kappa_{\hat Ric}$ and the diameter.
This constant depends continuously on its arguments  and 
increasingly on $\hat \kappa_{\hat Ric}$ and $diam_g(M)$.}\\

If we assume an upper bound $D$ for the diameter and a nonpositive  lower 
bound $\kappa$ for the Ricci 
curvature, then we have $\min_M \hat Ric  \ge D^2 \kappa$ and $C(n,p, \hat 
\kappa_{\hat Ric}, diam_g(M)) 
\le C(n,p, D^2  |\kappa|, D)$.  Hence the estimates (\ref{maxB}) can be  applied.    We state a corollary for the case of positive  
Ricci curvature.  We formulate it under the assumption $Ric \ge n-1$,  which 
can always 
be achieved by rescaling. \\

\noindent {\bf Theorem C} {\it Assume $n \ge 3$ and that the Ricci  
curvature satisfies $Ric \ge n-1$. 
Let  $u$ be a function in  $W^{1,\alpha}(M)$ with $\alpha>n$ satisfying 
\ba 
\Delta u \ge  f+\mbox{ div }\Phi
\ea
in the weak sense for a measurable function $f$ on  $M$ such that $f^- \in 
L^p(M)$  and a vector field
$\Phi \in L^{2p}(TM)$  with $p>\frac{n}{2}$.   Then we have 
\ba \label{maxC}
\sup_M  \, u \le u_M+ C(n, p)(\|f\|_p^*+\|\Phi\|^*_{2p})
\ea 
with a positive  constant $C(n, p)$ depending only on $n$ and $p$.}\\

Analogous results hold true if we assume an upper bound for the  diameter, 
and 
a lower bound for the Ricci curvature in a suitable integral  sense, thanks 
Gallot's  
and Petersen-Sprouse's estimates for the  Neumann isoperimetric constant in 
[Ga2] and 
[PS]. We omit the 
obvious  statements of those results. \\

\noindent {\bf Remark 1} In the above results we restrict to dimensions $n \ge 3$. 
The  2-dimensional analogues also hold true, see [Y]. We would also like to  
mention that it is straightforward to 
extend the above results to   compact manifolds with boundary under  
the Neumann boundary condition.  One can also extend the 
above results to general elliptic  operators of divergence form, see [Y]. \\

\noindent {\bf Remark 2}  Theorem A is also valid if we replace in the definition of 
$C^*_{I,N}(M,g)$ 
the 
Neumann isoperimetric constant $C_{N,I}(M,g)$ by the Poincar\'{e}-Sobolev 
constant  $C_{P,S}(M,g)$ (see Section  2 for its definition).  Indeed, 
it is 
the Poincar\'{e} inequality (\ref{poincare1}), the Poincar\'{e}-Sobolev 
inequality (\ref{poincare2}) and the Sobolev inequality (\ref{sobolev}) which are employed in 
our arguments. The Neumann isoperimetric constant appears in these inequalities. 
Obviously, the Poincar\'{e}-Sobolev inequality (\ref{poincare2}) can be reformulated in terms of 
the  Poincar\'{e}-Sobolev constant. Then the Poincar\'{e} inequality (\ref{poincare1}) and 
the Sobolev inequality (\ref{sobolev}) follow as corollaries, with the constants suitably modified. 
We formulate 
Theorem A in terms of the Neumann isoperimetric constant because we 
consider it to be a more fundamental quantity. \\

The proof of Theorem A  involves several 
ingredients. One is  
Moser iteration based on the Sobolev inequality. Various versions of this technique have been used in many  
situations, but 
the way it is done in this paper is new,  see 
the proof of 
Lemma \ref{maxth1}.  It is in this proof that the scaling invariance is broken, as mentioned above. 
On the other hand, from this  proof one can 
see that the  technique of Moser 
iteration alone  cannot  lead to a maximum estimate for $u-u_M$ in terms of 
$f$ and $\Phi$.  
Instead, the estimate one obtains also depends on the $L^2$ norm of  
$(u-u_M)^+$. 
Without using additional tools it seems impossible to go any  further.
Our strategy for overcoming this difficulty is to employ the Green  function 
$G_0$ of 
the Laplace operator. First we combine 
Lemma  \ref{maxth1} with the Poincar\'{e} inequality to establish Theorem 
\ref{maxth2}  which is 
the corresponding maximum principle for solutions (rather than  
subsolutions).  Using this result 
we reduce the right hand side of  (\ref{sub}) to a constant. 
Then we utilize the Green function $G_0$  to  obtain 
the desired estimate.  Employing the Green function is crucial  for the whole 
scheme.
 
There is an additional subtlety here. Usually, maximum principles based on  
Moser 
iteration hold true for all subsolutions $u$ in the Sobolev space  
$W^{1,2}(M)$.
(This is the case in Lemma \ref{maxth1} (for subsolutions) and  Theorem 
\ref{maxth2} (for solutions).)
In the  situation of Theorem A  (hence also Theorem B and Theorem C), we 
have to require $u \in  W^{1,\alpha}(M)$ for $\alpha>n$.
This restriction stems from the  involvement of the Green function. 
Using additional tools, one can extend  Theorem A to $u \in W^{1,2}(M)$, 
provided that 
$\Phi=0$, see [Y]. It remains  open whether one can extend the full Theorem A 
to 
$u \in W^{1,2}(M)$.   (See also [Y] for a weaker 
maximum principle which holds true for all $u \in  W^{1,2}(M)$.) 
 
In the above scheme of utilizing the Green function $G_0$, a lower  
bound for  $G_0$  is needed. In [Si], a lower 
bound for $G_0$  in terms of $\hat \kappa_{\hat Ric}$, the volume and the 
diameter is  obtained. This lower bound is sufficient for establishing the 
estimate  
(\ref{maxB}) in Theorem B and the estimate (\ref{maxC}) in Theorem C,  
but 
is not suitable for establishing the general estimate (\ref{maxA}) in  
Theorem A.  Following the  arguments in 
[CL] and [Si], we derive in Section 3  a lower bound for  $G_0$
which is proportional to  $vol_g(M)^{-1}$ with a factor given in terms of the volume-normalized Neumann isoperimetric constant. This form of lower  
bound is exactly 
what we need for establishing Theorem A. It is also of  
independent interest.

We would like to mention that Theorem C is sufficient for the purpose  of 
[P] 
because all involved functions in [P] are at least Lipschitz  continuous. We 
would also like to mention that  Theorem A
(or Theorem  \ref{maxth1} ) leads to an estimate for the $L^p$ norm of 
the Green function  $G_0$ for each $0<p<\frac{n}{n-2}$ (thanks an observation 
by Xiaodong  Wang) and an estimate for 
the $L^q$ norm of the gradient of $G_0$ for each  $0<q<\frac{n}{n-1}$. This 
will 
be presented elsewhere.

We would like to thank Xiaodong Wang for bringing the question  regarding 
the validity of Theorem C (with $\Phi=0$)  to 
our  attention. The first named author would also like to acknowledge many helpful discussions with Xiaodong Wang, and also with Jian Song.

\sect{The Neumann Isoperimetric Constant} 
 
Consider a closed  Riemannian manifold $(M, g)$ of dimension  $n$.
The Neumann isoperimetric constant of 
$(M, g)$  is defined to 
be 
\ba
C_{N,  I}(M,g)=\sup\{\frac{vol(\Omega)^{\frac{n-1}{n}}}{A(\partial \Omega)}: 
\Omega  \subset M \mbox{ is a } 
C^1 \mbox{ domain }, vol_g(\Omega) \le \frac{1}{2}  vol_g(M)\}.
\ea

The Poincar\'{e}-Sobolev constant (for the exponent $2$) of $(M,g)$ is defined 
to be 
\ba
C_{P,S}(M,g)=\sup\{\|u-u_M\|_{\frac{2n}{n-2}}: u \in C^1(M), \|\nabla u\|_2=1\}.
\ea 
 
We have the following Poincar\'{e} inequality,  Poincar\'{e}-Sobolev inequality 
and Sobolev inequality.  See [Y] for their proofs. (For 
these inequalities with somewhat different constants see [Si].)   
The Poincar\'{e}-Sobolev inequality (\ref{poincare2}) 
gives an upper bound of the Poincar\'{e}-Sobolev constant in terms of the 
Neumann isoperimetric constant.

\begin{lem} \label{poincare} There hold for all $u \in  W^{1,2}(M)$
\ba \label{poincare1}
\|u-u_M\|_2 \le  \frac{2(n-1)}{n-2}C_{N, I}(M, g)  
vol_g(M)^{\frac{1}{n}}\|\nabla  u\|_2,
\ea
\ba \label{poincare2}
\|u-u_M\|_{\frac{2n}{n-2}} \le  \frac{4(n-1)}{n-2} C_{N,I}(M,  g)
\|\nabla u\|_2,
\ea
and
\ba \label{sobolev}
\|u\|_{\frac{2n}{n-2}}  \le \frac{2(n-1)}{n-2} C_{N,I}(M,  g)\|\nabla  u\|_2+
\frac{\sqrt{2}}{vol_g(M)^{\frac{1}{n}}}  \|u\|_2,
\ea
whenever $n\ge  3$.
\end{lem}

It is convenient to use the following volume-normalized Neumann isoperimetric 
constant: 
\ba
C_{I,N}^*(M,g)=C_{N,I}(M,g)vol_g(M)^{\frac{1}{n}},
\ea
which was first introduced by J.~Cheeger in his study of the first eigenvalue of 
the Laplace operator [Che]. 
Note that $C_{I,N}(M,g)$ is scaling invariant, while $C^*_{I, M}(M,g)$ has the same scaling weight as 
the $n$-th root of the volume, or the 
diameter. In terms of $C^*_{I,N}(M,g)$ and the volume-normalized $L^p$ norms, Lemma \ref{poincare} can be reformulated as follows.

\begin{lem} \label{poincare*} There hold for all $u \in  W^{1,2}(M)$
\ba \label{poincare1*}
\|u-u_M\|^*_2 \le  \frac{2(n-1)}{n-2}C^*_{N, I}(M, g)\|\nabla  u\|^*_2,
\ea
\ba \label{poincare2*}
\|u-u_M\|^*_{\frac{2n}{n-2}} \le  \frac{4(n-1)}{n-2} C_{N,I}^*(M,  g)
\|\nabla u\|^*_2,
\ea
and
\ba \label{sobolev*}
\|u\|^*_{\frac{2n}{n-2}}  \le \frac{2(n-1)}{n-2} C_{N,I}^*(M,  g)\|\nabla  u\|^*_2+
\sqrt{2}  \|u\|^*_2,
\ea
whenever $n\ge  3$.
\end{lem}

The following estimate of the volume-normalized Neumann isoperimetric constant easily follows from  
S.~Gallot's corresponding estimate in [Ga1]. 
 
\begin{theo}  \label{gallot} There holds
\ba  \label{gallot}
C_{N,I}^*(g, M) \le C(n, \hat \kappa_{\hat  Ric})diam_g(M),
\ea
where $C(n, \hat \kappa_{\hat Ric})$ is a positive constant  depending only 
on $n$ and
$\hat \kappa_{\hat Ric}$.  It depends  continuously and increasingly on $\hat 
\kappa_{\hat Ric}$.
\end{theo} 
\Pf We rescale to make the diameter equal one.  Then we apply the 
estimate for the Neumann isoperimetric constant in [Ga1]. Expressing 
the estimate in terms of the original metric we arrive at (\ref{gallot}). 
\qed \\

\sect{The Green Function} 
 
Consider a closed Riemannian manifold $(M, g)$ of dimension $n$ as  before.  
Let $G_0(x, y)$ be the unique Green function of the Laplace  operator $\Delta$ 
such that $\int_M G_0(x, y)dy=0$ for all 
$x \in M$, where  $dy$ denotes the volume form of $g$.  Thus we have 
\ba  \label{greenrep}
u(x)=u_M-\int_M G_0(x, y)\Delta_y u(y)dy
\ea
for all  $u \in C^{\infty}(M)$, where $\Delta_y$ means $\Delta$ with the 
subscript  indicating 
the argument. (Similar notations will be used below.)   $G_0(x,y)$ is smooth 
away from 
$x=y$.  Moreover, $G_0(x,y)=G_0(y,x)$  
for all $x,y \in M, x \not =y$. In this section we present some basic  
facts about $G_0$ and derive a lower bound of $G_0$ in terms of 
$C_{N,  I}(M,g)$ and the volume.

\begin{lem} \label{greenintegral} Assume $n \ge 3$. Then there holds  
\ba
G_0(x,\cdot)  \in W^{1, \beta}(M)
\ea
for all  $0<\beta<\frac{n}{n-1}$.
\end{lem}
\Pf  By e.g. [Theorem 4.17,  A] we have $|G_0(x,y)| \le Cd(x,y)^{2-n}$ and 
$|\nabla_y G_0(x,y)| \le  Cd(x,y)^{1-n}$ for all $x, y \in M, x \not = y$ and 
a positive constant $C$  depending on $g$. Consequently, we have  
\ba  \label{greenintegral-1}
G_0(x, \cdot) \in L^{q_1}(M), \, \nabla_y G_0(x,y)  \in L^{q_2}(M)
\ea 
for all $0<q_1<\frac{n}{n-2}$ and   all
$0<q_2<\frac{n}{n-1}$. Then it follows easily that 
$G(x, \cdot)  \in W^{1, p}(M)$ for all $x \in M$ and 
$0<q<\frac{n}{n-1}$.   Indeed, we have for an arbitary $x \in M$ and small 
$\epsilon>0$  
\ba
\int_{M-B_{\epsilon}(x)} G_0(x,y) \mbox{ div}_y \Phi(y)  dy
&=&-\int_{M-B_{\epsilon}(x)} \nabla_y G_0(x,y) \cdot \Phi(y)dy  \nonumber \\
&+&\int_{\partial B_{\epsilon}(x)} 
G_0(x,y) \Phi(y)  \cdot \nu(y)
\ea
for all smooth vector fields $\Phi$ on $M$, where $\nu$  denotes the inward 
unit normal of the geodesic 
sphere
$\partial  B_{\epsilon}(x)$. Since $|G_0(x, y)| \le C\epsilon^{2-n}$ on 
$\partial  B_{\epsilon}(x)$ we can let  $\epsilon \rightarrow 0$ to arrive at 
\ba  \label{greenweak}
\int_{M} G_0(x,y) \mbox{ div}_y \Phi(y) dy
=-\int_{M}  \nabla_y G_0(x,y) \cdot \Phi(y)dy.
\ea
By (\ref{greenintegral-1}) and  (\ref{greenweak}) we infer that 
$G(x, \cdot) \in W^{1,q}(M)$ for all  $0<q<\frac{n}{n-1}(M)$ and 
that (\ref{greenweak}) holds true for all  $\Phi \in W^{1, p}(TM)$ 
whenever $p>n$, where $W^{1,p}(TM)$ denotes the  $W^{1,p}$ Sobolev 
space of vector fields on $M$.  \\ \qed \\

\begin{lem} \label{greennew} Let $u \in W^{1,q}(M)$ with $q >n$.  Then 
\ba \label{greenrep1}
u(x)=u_M+\int_M \nabla_y G_0(x, y) \cdot  \nabla_y u (y) dy
\ea
holds true for a.e. $x \in M$. 
\end{lem}
\Pf  By Lemma \ref{greenintegral}, we 
can integrate (\ref{greenrep}) by parts to  deduce (\ref{greenrep1}) for 
all $u \in C^{\infty}(M)$. Applying Lemma  \ref{greenintegral} and 
a limiting argument
we then conclude that  (\ref{greenrep1}) holds true for each 
$u \in W^{1,q}(M)$ a.e. as long as  $q>n$. \\ \qed \\
 
Next let $H(x,y,t)$ be the heat kernel for $\Delta$, i.e. 
\ba  \label{kernelequation}
\frac{\partial}{\partial t}H(x,y, t)=\Delta_y H(x, y,  t)
\ea
for $t>0$ and 
\ba \label{delta} \lim_{t \rightarrow 0} H(x,  y, t)=\delta_x
\ea
in the sense of distributions, where $\delta_x$ is the  Dirac 
$\delta$-function with center $x$. 
$H$ is symmetric in $x, y$  and  smooth away from $x=y, t=0$.   We have the 
basic  
representation formula
\ba   \label{heatrep}
u(x,t)=\int_{0}^td\tau \int_M H(x,y,t-\tau)  (\frac{\partial}{\partial 
\tau}-\Delta_y) u(y, \tau) dy+\int_M H(x,y,t) 
u(y,  0)dy
\ea
for all smooth functions $u$ and $t>0$. 
Note that   
$H(x,y,t)>0$ for $t>0$ and all $x, y \in M$. We  set
\ba
G(x,y,t)=H(x,y,t)-\frac{1}{vol_g(M)}.
\ea
Choosing  $u(x,t)\equiv 1$ in (\ref{heatrep}) we deduce 
\ba \label{H-1}
\int_M  H(x,y,t)dy =1
\ea
and hence 
\ba \label{G-0}
\int_M  G(x,y,t)dy=0
\ea
for all $x \in M$ and $t>0$.

\begin{lem} Assume $n \ge 3$. Then there holds 
\ba  \label{greenequalA}
G_0(x,y)=\int_0^{\infty}  G(x,y,t)dt.
\ea
\end{lem}
\Pf We have 
\ba  \label{infty}
|H(x,y,t)-\frac{1}{vol_g(M)}| \le  Ct^{-\frac{n}{2}}
\ea
for a certain positive constant $C$ depending on $g$  (for a geometric 
estimate of $C$ see [CL]). 
On the other hand, we have   the inequality (see e.g. [Proposition VII.3.5, 
Ch])
\ba  \label{zero}
|\frac{H(x,y,t)}{{\cal H}(x,y,t)}-1| \le Cd(x,y) 
\ea
for  all $t>0$ and $x, y \in M$ with $d(x,y) \le \frac{1}{4}inj_g(M)$, where 
$C$  is a positive 
constant depending on $g$,  $inj_g(M)$ denotes the  injectivity radius of 
$(M,g)$, and 
\ba
{\cal H}(x,y,t)=(4\pi  t)^{-\frac{n}{2}} e^{-\frac{d(x,y)^2}{4t}}.
\ea

By (\ref{heatrep}) we have for a smooth function $u(x)$
\ba  \label{heatrep1}
u(x)&=& -\int_0^{t} d\tau\int_M H(x,y,t-\tau)  \Delta_y u(y) dy+\int_M 
H(x,y,t)u(y)dy \nonumber \\
&=& -\int_0^{t}  d\tau\int_M G(x,y,t-\tau) \Delta_y u(y) dy +\int_M 
G(x,y,t)u(y)dy +u_M \nonumber  \\ 
&=& -\int_0^{t} ds\int_M G(x,y,s) \Delta_y u(y) dy +\int_M  G(x,y,t)u(y)dy 
+u_M 
\ea
By (\ref{infty}) and (\ref{zero}) we can  let $t\rightarrow \infty$ in 
(\ref{heatrep1}) to  obtain
\ba
u(x)=-\int_M (\int_0^{\infty} G(x,y,s)ds) \Delta_y u(y)  dy
+u_M.
\ea
On the other hand, by (\ref{G-0}) we deduce  
\ba
\int_M (\int_0^{\infty} G(x,y,s)ds)dy=0
\ea
for all $x \in  M$.  We conclude that (\ref{greenequalA}) holds true. \\ 
\qed \\

\begin{lem} \label{greenformula} There holds 
\ba  \label{greenequal}
G(x,y,t+s)=\int_M G(x,z,s)G(z,y,t)dz
\ea
for all $x,  y \in M$ and $t>0, s>0$, where $dz$ denotes the volume form of 
$g$ with $z  \in M$ as 
the argument. In particular, we have 
\ba  \label{squareformula}
G(x,x,t)= \int_M G(x,y,\frac{t}{2})^2dy 
\ea
and  it follows that $G(x,x,t)>0$ for all $x \in M$ and $t>0$. 
\end{lem}  
\Pf   Note that $G(x,y,t)$ satisfies the heat equation 
\ba  \label{G-heat}
\frac{\partial}{\partial t}G(x,y,t)=\Delta_y G(x,y,t)=\Delta_x  G(x,y,t).
\ea
Choosing $u(x,t)=G(x,y,t+s)$  in (\ref{heatrep}) for each fixed $y$  we  deduce, on account 
of 
(\ref{G-heat}) and (\ref{G-0})  
\ba
G(x,y,t+s)= \int_M H(x,z, t)G(z,y,s)dz= \int_M  G(x,z,t)G(z,y,s)dz.
\ea
Switching $t$ with $s$ we arrive at the desired equation (\ref{greenequal}).
 
The formula (\ref{squareformula}) follows immediately and hence $G(x,x,t)  
\ge 0$. If 
$G(x,x,t_0)=0$ for some $x$ and $t_0>0$, then  (\ref{squareformula}) implies 
that 
$G(x,y,\frac{t_0}{2})=0$ for all $y \in  M$. Then $G(x,y, t)=0$ for all $y 
\in M$ 
and $t \ge \frac{t_0}{2}$, because  $G(x,y,t)$ satisfies the 
heat equation.
It follow that  $H(x,y,t)=vol_g(M)^{-1}$ for all $y \in M$ 
and $t \ge \frac{t_0}{2}$. This  contradicts (\ref{heatrep}) as is easy to 
see. We conclude that  
$G(x,x,t)>0$ for all $x\in M$ and $t>0$.
\\ \qed \\

\begin{theo} \label{greenth} Assume $n \ge 3$. Then there holds 
\ba  \label{greenestimate}
G_0(x,y) \ge -C_0(n) C_{I,N}^*(M,g)^2 vol_g(M)^{-1}
\ea
for all $x, y \in M, x \not =y$, where  
\ba
C_0(n)=\frac{8n^2(n-1)^2}{(n-2)^3}
\left(\frac{n-2}{2}\right)^{\frac{4}{n}}.
\ea
\end{theo}
\Pf This  follows from the arguments in [CL] and [Si] with some 
modification. By the  rescaling invariance of (\ref{greenestimate}) 
we can assume   $vol_g(M)=1.$
Differentiating the equation (\ref{greenequal}), 
setting $y=x$  and  
replacing $t$ and $s$ by $\frac{t}{2}$ we deduce for  $t>0$
\ba \label{green-t}
\frac{\partial G}{\partial t}(x,x,t)=\int_M  \frac{\partial G}{\partial 
t}(x,z, \frac{t}{2}) 
G(x,z, \frac{t}{2})  dz=\int_M (\Delta_z G(x, z, \frac{t}{2})) G(x, z, 
\frac{t}{2}) dz.
\ea
We  integrate  (\ref{green-t}) by parts to derive 
\ba  \label{green-nabla}
\frac{\partial G}{\partial t}(x,x,t)=-\int_M   |\nabla_z G(x, z, 
\frac{t}{2})|^2  dz.
\ea
 
Applying the Poincar\'{e}-Sobolev inequality (\ref{poincare2}) we then obtain  
\ba
-\frac{\partial G}{\partial t}(x,x,t) \ge  
\left(\frac{4(n-1)}{n-2}C_{N,I}(M, g)\right)^{-2} 
\left(\int_M  |G(x,z,\frac{t}{2})|^{\frac{2n}{n-2}}dz \right)^{\frac{n-2}{n}}.
\ea 
By  H\"{o}lder's inequality we have
\ba
\left(\int_M |G(x,z,\frac{t}{2})|^2dz  \right)^{\frac{n+2}{n}} &=& 
\left(\int_M  |G(x,z,\frac{t}{2})|^{\frac{4}{n+2}}|G(x,z, 
\frac{t}{2})|^{\frac{2n}{n+2}}  dz
\right)^{\frac{n+2}{n}} \nonumber \\
&\le&  
\left(\int_M  |G(x,z,\frac{t}{2})|^{\frac{2n}{n-2}}dz\right)^{\frac{n-2}{n}}
\left(\int_M  |G(x,z,\frac{t}{2})|dz\right)^{\frac{4}{n}}. \nonumber \\
\ea
Next observe  that $\int_M |G(x, z,t)|dz \le 2$ because $H(x, z, t)>0$.
Hence we arrive  at 
\ba \label{greensquare}
-\frac{\partial G}{\partial t}(x,x,t)  \ge  C
\left(\int_M |G(x,z,\frac{t}{2})|^{2}dz \right)^{\frac{n+2}{n}}=C  
G(x,x,t)^{\frac{n+2}{n}},
\ea
where 
\ba
C=\left(  \frac{4(n-1)}{n-2} C_{N,I}(M,g)\right)^{-2}.
\ea
Integrating  (\ref{greensquare}) we derive
\ba
G(x,x,t)^{-\frac{2}{n}} \ge  G(x,x,s)^{-\frac{2}{n}}+\frac{n}{2}C (t-s)
\ea
for $t>s>0$.   (Note that $G(x,x,t)>0$ by Lemma \ref{greenformula}.)
Letting $s  \rightarrow 0$  we infer $G(x,x,t)^{-\frac{2}{n}} \ge 
\frac{n}{2}Ct$, and  hence
$G(x,x,t) \le C^{-\frac{n}{2}}(\frac{n}{2})^{\frac{n}{2}}  t^{-\frac{n}{2}}$. 
Now we have by Lemma  \ref{greenformula}
\ba
|G(x,y,t)| &=& |\int_M  G(x,z,\frac{t}{2})G(z,y,\frac{t}{2})dz| \nonumber  \\
&\le&
\left(\int_M G(x,z,\frac{t}{2})^2dz  \right)^{\frac{1}{2}} 
\left(\int_M G(z,y,  \frac{t}{2})^2dz\right)^{\frac{1}{2}} \nonumber \\
&=&  G(x,x,\frac{t}{2})^{\frac{1}{2}} 
G(y,y,\frac{t}{2})^{\frac{1}{2}} \le  
C^{-\frac{n}{2}}(\frac{n}{2})^{\frac{n}{2}} t^{-\frac{n}{2}}.
\ea
 
Since $H(x,y,t)>0$ we have $G(x,y,t) \ge -\frac{1}{vol_g(M)}=-1$. We  deduce
for each $\tau>0$ 
\ba \label{G}
G(x,y)&=&  \int_0^{\infty}G(x,y,t)dt \ge -\int_0^{\tau} dt 
-\int_{\tau}^{\infty}  C^{-\frac{n}{2}}(\frac{n}{2})^{\frac{n}{2}} t^{-\frac{n}{2}}
\nonumber  \\
&=& -\tau-\frac{n-2}{2}C^{-\frac{n}{2}}(\frac{n}{2})^{\frac{n}{2}}  
\tau^{-\frac{n-2}{2}}=-\tau-C_1\tau^{-\frac{n-2}{2}},
\ea
where $C_1=\frac{n-2}{2}C^{-\frac{n}{2}}(\frac{n}{2})^{\frac{n}{2}}$.
The maximum of the function $\tau+C_1\tau^{-\frac{n-2}{2}}$ is achieved at 
$\tau=({C_1(n-2)/2})^{2/n}$ and hence equals  
\ba
\frac{n}{n-2} (C_1\frac{n-2}{2})^{\frac{2}{n}}=\frac{8n^2(n-1)^2}{(n-2)^3}
\left(\frac{n-2}{2}\right)^{\frac{4}{n}}C_{I, N}(M,g)^2.\nonumber
\ea
we arive at 
\ba
G(x,y)  \ge -\frac{8n^2(n-1)^2}{(n-2)^3}
\left(\frac{n-2}{2}\right)^{\frac{4}{n}}C_{I, N}(M,g)^2,
\ea
which leads  to (\ref{greenestimate}) by rescaling. \\

\sect{Neumann Type Maximum Principles}

In this section we consider a fixed closed Riemannian manifold $(M, g)$ of dimension 
$n \ge 3$ as before. 

\begin{lem} \label{holder} 1) There hold for $f_1 \in L^p(M)$ and $ f_2 \in L^q(M)$ with 
$p^{-1}+q^{-1}=1$
\ba \label{holder1}
\|f_1f_2 \|_1^* \le \|f_1\|^*_p \cdot \|f_2\|^*_q.
\ea
2) There holds for $p \ge q \ge 1$ and $f \in L^p(M)$
\ba \label{holder2}
\|f\|_q^* \le \|f\|_p^*.
\ea
\end{lem}
\Pf These follow straightforwardly from the classical H\"{o}lder inequality. \qed \\

\begin{lem} \label{maxth1} Let $n \ge 3$.  Assume that  $u \in  W^{1,2}(M)$ 
satisfies 
\ba \label{delta}
\Delta u \ge f+\mbox{ div }  \Phi
\ea
in the weak sense, i.e. 
\ba \label{nabla}
\int_M \nabla u  \cdot \nabla \phi \le  
-\int_M f \phi + \int_M \Phi \cdot \nabla  \phi
\ea
for all {\it nonnegative} $\phi \in W^{1,2}(M)$, where $f$ is a  measurable 
function on $M$ 
such that $f^- \in L^p(M)$, and $\Phi \in  L^{2p}(TM)$, with $p>\frac{n}{2}$. 
   Then we have  
\ba  \label{max1}
\sup_M \, (u-\lambda) &\le&
A C_1  (\|f^-\|^*_p
+\|\Phi\|_{2p}^*) \nonumber  \\
&&+A(C_1+\sqrt{2}) \|(u-\lambda)^+\|^*_2
\ea
for each  $\lambda \in {\bf R}$, where $A$ and $C_1$ are positive numbers 
depending only  on $n, p$ and $C_{N,I}^*(M,g)$. 
Their explicit values are given in the proof  below.
\end{lem}
\Pf The arguments here are inspired by some arguments in  [GT]. A special new feature in our 
argument is the use of the volume-normalized isoperimetric constant (or Sobolev 
constant) and the volume-normalized $L^p$ norms.  We set 
\ba
a=\|f^-\|^*_p+\|\Phi\|^*_{2p}.
\ea
Then we  set $b=a$ if $a>0$ and $b=1$ if $a=0$. 
For $L>|\lambda|$  we set  $w=\min\{(u-\lambda)^+, L\}+b$.  Then $w \in 
W^{1,2}(M)$ and is bounded. It  follows that  
$w^{\gamma} \in W^{1,2}(M)$ for $\gamma \ge 1$.   Moreover, we have 
$\nabla w=0$ if $u \ge \lambda+ L$, $\nabla w=\nabla u$ if  
$\lambda<u<L+\lambda$, and $\nabla w=0$ if $u \le \lambda$.    Choosing 
$\phi=w^{\gamma}({\vol M})^{-1}$ with ${\vol M}=vol_g(M)$ in (\ref{nabla}) we obtain
\ba
\frac{ \gamma}{\vol M} \int_M  |\nabla w|^2 w^{\gamma-1} &\le& -\frac{1}{\vol M}\int_M f w^{\gamma} + \frac{ \gamma }{\vol M} 
\int_M  w^{\gamma-1}\Phi \cdot \nabla w \nonumber \\
 &\le& \frac{1}{\vol M} \int_M |f^-|  w^{\gamma} + \frac{ \gamma }{\vol M}
\int_M w^{\gamma-1}\Phi \cdot \nabla w.
\ea
 
First we choose $\gamma=1$ to deduce
\ba \label{gamma=1}
\|\nabla w\|_2^{*2}  &\le& \| f^- w \|_1^* +\|  \Phi \cdot \nabla w\|_1^*  \nonumber 
 \\
&\le& \frac{1}{2}\|f^-\|_2^{*2}+ \frac{1}{2}\|w\|_2^{*2} +  \frac{1}{2}\|\Phi\|_2^{*2} 
+ \frac{1}{2} \|\nabla w\|_2^{*2},
\ea
where we have used Lemma~\ref{holder}. 
It follows  that
\ba \label{gamma=11}
\|\nabla w\|_2^{*2} \le \|f^-\|_2^{*2}+\|\Phi\|_{2}^{*2}+ \|w\|_2^{*2} \le  
\|f^-\|_p^{*2}+\|\Phi\|_{2p}^{*2}+ \|w\|_2^{*2},
\ea
where the second inequality  follows from Lemma \ref{holder}.
Applying  the Sobolev inequality (\ref{sobolev*}) we then deduce 
\ba  \label{w-2n/n-2}
\|w\|^*_{\frac{2n}{n-2}}  &\le&  \frac{2(n-1)}{n-2}C^*_{N,I}(M,g)
\|\nabla  w\|^*_2+\sqrt{2}\|w\|^*_2   \nonumber \\
&\le&  C_1  (\|f^-\|^*_p+\|\Phi\|^*_{2p})
+(C_1+\sqrt{2}) \|w\|^*_2,
\ea
where  $C_1=\frac{2(n-1)}{n-2}C^*_{N,I}(M,g)$.

Next we consider general $\gamma\ge 1$. We deduce 
\ba
&&\frac{\gamma}{\vol M} \int_M |\nabla w|^2 w^{\gamma-1} \le  \frac{1}{b \vol (M)}\int_M |f^-|  
w^{\gamma+1} 
+\frac{\gamma}{b \vol (M)}\int_M w^{\gamma} |\Phi| |\nabla w|  \nonumber  \\
&&\le \frac{1}{b} \|f^- w^{\gamma+1}\|_1^*  +\frac{\gamma}{2}
\| |\nabla w|^2 w^{\gamma-1} \|_1^*+ \frac{\gamma}{2b^2}  
\| w^{\gamma+1} |\Phi|^2 \|_1^*.
\ea
It follows that, on account of Lemma \ref{holder}
\ba
&&  \frac{\gamma}{2}\| |\nabla w|^2 w^{\gamma-1} \|_1^*  \nonumber  \\
&\le&  \frac{1}{b} \| |f^- w^{\gamma+1}\|_1^*
+  \frac{\gamma}{2b^2} 
\| w^{\gamma+1} |\Phi|^2 \|_1^* \nonumber  \\
&\le& \frac{1}{b}\|f^-\|_p^* \cdot  \|w\|_{(\gamma+1)\frac{p}{p-1}}^{*(\gamma+1)} 
+  \frac{\gamma}{2b^2}    {\|w\|}_{(\gamma+1)\frac{p}{p-1}}^{*(\gamma+1)} \cdot
\|\Phi\|^{*2}_{2p} \nonumber \\
&\le&   \|w\|_{(\gamma+1)\frac{p}{p-1}}^{*(\gamma+1)} 
+  \frac{\gamma}{2} \|w\|_{(\gamma+1)\frac{p}{p-1}}^{*(\gamma+1)}.
\ea
It follows that
\ba \label{moser1}
\|\nabla  w^{\frac{\gamma+1}{2}}\|_2^{*2} 
\le  
\frac{(\gamma+2)(\gamma+1)^2}{4\gamma}\|w\|_{(\gamma+1)\frac{p}{p-1}}^{*(\gamma+1)}.
\ea
 
Now we apply the Sobolev inequality (\ref{sobolev*}) and (\ref{moser1}) to deduce
\ba  \label{moser2}
\|w\|_{(\gamma+1) \frac{n}{n-2}}^{*\frac{\gamma+1}{2}}  
&=& \|w^{\frac{\gamma+1}{2}}\|^*_{\frac{2n}{n-2}}\le  
A_{\gamma}\|w\|_{(\gamma+1)\frac{p}{p-1}}^{*\frac{\gamma+1}{2}}
+\sqrt{2}\|w^{\frac{\gamma+1}{2}}\|^*_2  \nonumber \\
&=&  
A_{\gamma}\|w\|_{(\gamma+1)\frac{p}{p-1}}^{*\frac{\gamma+1}{2}}+\sqrt{2}\|w\|_{\gamma+1}^{*\frac{\gamma+1}{2}}
  \nonumber \\
&\le&  
(A_{\gamma}+\sqrt{2})\|w\|_{(\gamma+1)\frac{p}{p-1}}^{*\frac{\gamma+1}{2}},
\ea
where  
\ba
A_{\gamma}=C^*_{N,I}(M,g)
\frac{(n-1)(\gamma+1)}{n-2}  \sqrt{\frac{\gamma+2}{\gamma}}.
\ea
Consequently, we obtain  
\ba
\|w\|^*_{(\gamma+1) \frac{n}{n-2}} \le  
(A_{\gamma}+\sqrt{2})^{\frac{2}{\gamma+1}}\|w\|^*_{(\gamma+1)\frac{p}{p-1}}.
\ea
Replacing  $\gamma+1$ by $\gamma \ge 2$ we infer 
\ba \label{first}
\|w\|^*_{\gamma  \frac{n}{n-2}} \le  
(A_{\gamma-1}+\sqrt{2})^{\frac{2}{\gamma}}\|w\|^*_{\gamma\frac{p}{p-1}}.
\ea

Now we choose $ \gamma_0=1+\frac{n(p-2)+2p}{(n-2)p}$ and  
$\gamma_k=\gamma_{k-1} \frac{n(p-1)}{(n-2)p}$ for 
$k \ge 1$, i.e.  $\gamma_k=\gamma_0 (\frac{n(p-1)}{(n-2)p})^k$.  Since 
$p>\frac{n}{2}$,  we have $\gamma_0>2$ and $\frac{n(p-1)}{(n-2)p}>1$. 
We also have $  \gamma_0 \frac{p}{p-1}=\frac{2n}{n-2}$.  We deduce
\ba  \label{second}
\|w\|^*_{\gamma_k} \le \left(\prod\limits_{1\le i\le k}  
(A_{\gamma_i-1}+\sqrt{2})^{\frac{2}{\gamma_i}}\right)  \|w\|^*_{\frac{2n}{n-2}}.
\ea
Since $\frac{n(p-1)}{(n-2)p}>1$, the  product $\prod\limits_{1\le i<\infty}  
(A_{\gamma_i-1}+\sqrt{2})^{\frac{2}{\gamma_i}}$
converges. We denote its  value by $A$. Letting $k \rightarrow \infty$ we 
infer, on account 
of  (\ref{w-2n/n-2})
\ba
\|w\|_{\infty} \le A\|w\|^*_{\frac{2n}{n-2}} \le A C_1  
(\|f^-\|^*_p+\|\Phi\|^*_{2p})
+A(C_1+\sqrt{2}) \|w\|^*_2.
\ea
Letting $L  \rightarrow \infty$ we then arrive at (\ref{max1}). \\
\qed \\

\noindent {\bf Remark 3} An important point in the above proof is to 
break the scaling invariance. Basically, the construction of the 
function $w$ is not scaling invariant. More precisely, if $g$ is transformed to 
 $\bar g=\alpha g$ for 
a positive constant $\alpha$, then (\ref{delta}) is tranformed to 
\ba
\Delta_{\bar g} u \ge \bar f+div \bar \Phi,
\ea
where $\bar f=\alpha^{-1}f$ and $\bar \Phi=\alpha^{-1} \Phi$.  We have
\ba
\|\bar f\|^*_{p, \bar g}+ \|\bar \Phi\|^*_{2p, \bar g}=\alpha^{-1} (
\|\bar f\|^*_{p}+ \|\bar \Phi\|^*_{2p}).
\ea
It follows that $a$ and hence $b$ are not scaling invariant. 
The fact that the estimate (\ref{max1}) is 
not scaling invariant is a result of this. This unconventional feature is needed for our purpose of 
controlling the constants in the estimates in terms of $C_{I, N}^*(M, g)$ alone. \\ 

\noindent {\bf Remark 4} Since the estimate (\ref{max1}) is not scaling invariant, one may wonder 
what happens to it if one lets the above scaling factor $\alpha$ go to $0$ or $\infty$. 
The answer to this question is simple:
the estimate deteriorates in the process. Indeed, as $\alpha \rightarrow 0$, the factor 
$AC_1$ in the estimate converges to a positive constant depending only on the dimension, but 
(the transformed) $\|f\|_p^*+\|\Phi\|^*_{2p}$ approaches $\infty$. As $\alpha \rightarrow 
\infty$, $AC_1$ approaches $\infty$ and $\|f\|_p^*+\|\Phi\|^*_{2p}$ approaches $0$, but 
the former has more weight than the latter.   
On the other hand, the non-invariance of the estimate allows us to vary $\alpha$ to obtain the 
optimal estimate. We do not pursue this in this paper because it is not needed for our 
main purpose. 

The same question can be asked about the estimate in Theorem A.  The answer is obviously the same. \\

\begin{theo} \label{maxth2} Let $n \ge 3$. Assume that 
$u \in  W^{1,2}(M)$, $f \in L^p(M)$ and $\Phi \in L^{2p}(TM)$ for some  
$p>\frac{n}{2}$, which satisfy
\ba \label{delta2}
\Delta u=f+\mbox{ div  } \Phi
\ea
in the weak sense, i.e. 
\ba \label{nabla2}
\int_M \nabla  u \cdot \nabla \phi = 
-\int_M f \phi + \int_M \Phi \cdot \nabla  \phi
\ea
for all $\phi \in W^{1,2}(M)$.    Then we have  
\ba  \label{max2}
\sup_M \, |u-u_M| \le C_2 (\|f\|^*_p+\|\Phi\|^*_{2p}),
\ea
where  $C_2=AC_1\left[1+2\max\{C_1,1\}(C_1+\sqrt{2})\right]$ with
$A$ and $C_1$  being from Lemma \ref{maxth1}.
\end{theo} 
\Pf 
Choosing $\phi=(u-u_M)({\vol M})^{-1}$ in  (\ref{nabla2}) we deduce 
by applying Lemma \ref{holder}, as in  (\ref{gamma=1}) and (\ref{gamma=11})
\ba 
\frac{1}{\vol(M)} \int_M |\nabla u|^2 &=& -\frac{1}{\vol(M)}\int_M f(u-u_M)+ \frac{1}{\vol(M)} \int_M  \Phi \cdot (u-u_M) 
\nonumber \\
&\le&   \|f\|^*_2 \cdot \|u-u_M\|^*_2 
+\|\Phi\|^*_{2} \cdot \|\nabla u\|^*_2 \nonumber \\  
&\le&  \|f\|^*_p \cdot \|u-u_M\|_2 +\|\Phi\|^*_{2p} \cdot \|\nabla u\|_2  \nonumber \\
&\le& \|f\|^*_p \cdot \|u-u_M\|^*_2 +\frac{1}{2} \|\Phi\|^{*2}_{2p}  +\frac{1}{2} 
\|\nabla u\|_2^{*2}.
\ea  
Hence 
\ba
\|\nabla u\|_2^{*2}  \le 2\|f\|^*_p \cdot \|u-u_M\|^*_2 +\|\Phi\|^{*2}_{2p}. 
\ea
Combining this  with the Poincar\'{e} inequality (\ref{poincare1*})  we then obtain
\ba
\|u-u_M\|^*_2 &\le&  
\frac{2(n-1)}{n-2}C^*_{N,I}(M,g) \left(\sqrt{2}\|f\|^{*\frac{1}{2}}_p \cdot
\|u-u_M\|^{*\frac{1}{2}}_2+\|\Phi\|^*_{2p}\right)
\nonumber \\ &\le&  \frac{1}{2} \|u-u_M\|^*_2 +C_1^2\|f\|^*_p
+C_1\|\Phi\|^*_{2p},
\ea
where  $C_1=\frac{2(n-1)}{n-2}C^*_{N,I}(M,g)$ as before.
It follows that 
\ba  \label{square}
\|u-u_M\|^*_2 \le 2C_1^2\|f\|^*_p 
+2C_1\|\Phi\|^*_{2p} \le  2\max\{C_1^2,C_1\}(\|f\|^*_p+\|\Phi\|^*_{2p}).
\ea
Combining (\ref{max1}) with  $\lambda=u_M$ and (\ref{square}) we then arrive 
at 
\ba
\sup_M \, (u-u_M)  \le C_2 (\|f\|^*_p+ \|\Phi\|^*_{2p}).
\ea
Replacing $u$ by $-u$  we obtain 
\ba
\inf_M \, (u-u_M) \ge -C_2 (\|f\|^*_p+  \|\Phi\|^*_{2p}).
\ea
The estimate (\ref{max2}) follows. \\ \qed \\

\noindent {\bf Proof of Theorem A} \\
 
Replacing $f$ by $f^-$ we can 
assume $f\le 0$. There is  a unique weak solution $v \in W^{1, 2p}(M)$ of the 
equation
\ba  \label{help}
\Delta v=f-f_M+\mbox{ div } \Phi
\ea
with $v_M=0$. Indeed,  we can minimize the functional
\ba \label{F}
F(v)=\int_M (|\nabla  v|^2-(f-f_M)v-\Phi \cdot \nabla v) 
\ea
for $v \in W^{1,2}(M)$ under  the constraint $v_M=0$.  By the 
H\"older inequality and the  Poincar\'{e} inequality (\ref{poincare1})
we have
\ba
F(v) \ge  c\|v\|_{1,2}^2-C(\|f-f_M\|_p^2 +\|\Phi\|_{2p}^2)
\ea
for some positive  constants $c$ and $C$, where $\|v\|_{1,2}$ denotes 
the $W^{1,2}$ norm of  $v$.  Hence a minimizer $v$ exists, which
is a desired solution of  (\ref{help}).  Its uniqueness follows from Theorem 
\ref{maxth2}.
The  property $v \in W^{1,2p}(M)$ follows from the regularity theory for 
elliptic  operators in divergence form.   
By Theorem \ref{maxth2} we have  
\ba \label{v}
\sup_M \, v \le C_2 (\|f-f_M\|^*_p + \|\Phi\|^*_{2p}) \le  C_2(2\|f\|^*_p+\|\Phi\|^*_{2p}) .
\ea
We set $w=u-v$. Then we have $w \in W^{1, q}(M)$ with  
$q =\min\{2p, \alpha\}>n$. There holds
\ba
\Delta w =\Delta u  -\Delta v \ge f_M
\ea
in the weak sense, i.e.
\ba  \label{weak}
\int_M \nabla w \cdot \nabla \phi  \le -f_M \int_M  \phi
\ea
for all nonnegative $\phi \in W^{1, \frac{q}{q-1}}(M)$. 
Now  we apply (\ref{greenrep1}) 
to $w$ to deduce
\ba
w(x)=w_M+ \int_M  \nabla_y G_0(x, y) \cdot \nabla_y w(y) dy
\ea
for a.e. $x \in M$.
Set $\sigma=\inf_{x \not =y} G_0(x, y)$. By 
(\ref{greenestimate}) there holds 
\ba
\sigma \ge -C_0(n) C^*_{I,N}(M,g)^2 vol_g(M)^{-1}.
\ea
Next we set
\ba
G (x,y)=G_0(x,y)-\sigma.
\ea
We have $G(x,y) \ge 0$ and 
\ba
w(x)=w_M+  \int_M \nabla_y G(x, y) \cdot \nabla_y w(y) dy
\ea
for a.e. $x \in  M$.   
By  Lemma \ref{greenintegral} and the fact   $\frac{q}{q-1} <\frac{n}{n-1}$, 
(\ref{weak}) holds true with  
$\phi=G(x, \cdot)$ for each given $x\in M$, i.e. 
\ba 
\int_M \nabla_y  G(x, y) \cdot \nabla_y w dy \le -f_M \int_M G(x, y)dy.
\ea
We then  deduce
\ba \label{w}
w(x)-w_M &\le& -f_M \int_M G(x,y)dy \nonumber  \\
&=& |f_M| (\int_M G_0(x,y)dy) -\sigma vol_g(M)) \nonumber \\  &=&
-\sigma vol_g(M) |f_M| \nonumber \\ &\le& C_0(n) C^*_{I,N}(M,g)^2 |f_M| 
\le C_0(n) C^*_{I,N}(M,g)^2 \|f\|^*_1 
\nonumber \\ &\le& C_0(n) C^*_{I,N}(M,g)^2  \|f\|^*_p
\ea
for  a.e. $x \in M$. 
Since $w_M=u_M$, combining (\ref{v}) and (\ref{w}) yields  
\ba
u(x)\le u_M+(C_0(n)C^*_{I,N}(M,g)+2C_2) \|f\|^*_p + C_2\|\Phi\|^*_{2p}
\ea
for a.e. $x \in M$. We  arrive at (\ref{maxA}) (note that 
$\sup_M u$ means the essential supremum).  
\qed \\
 
\noindent {\bf Proof of Theorem B} \\
 
The estimate (\ref{maxB}) follows straightforwardly from Theorem A and Theorem \ref{gallot}.  \\
\qed \\

\noindent {\bf Proof of Theorem C} \\
 
We have $\hat \kappa_{\hat Ric}=0$. By  Bonnet-Myers Theorem we have 
$diam_g(M) \le \pi$.  Hence Theorem B implies 
\ba  
u(x) \le u_M+ C(n, p, 0, \pi) (\|f^-\|_p^*+
\|\Phi\|_{2p}^*). \nonumber 
\ea 
  \qed \\

\end{document}